# GENERALIZATIONS OF CEVA'S THEOREM AND APPLICATIONS

Florentin Smarandache
University of New Mexico
200 College Road
Gallup, NM 87301, USA
E-mail: smarand@unm.edu

In these paragraphs one presents three generalizations of the famous theorem of Ceva, which states:

"If in a triangle $ABC$ one draws the concurrent straight lines $AA_1$, $BB_1$, $CC_1$, then $\dfrac{\overline{A_1B}}{\overline{A_1C}} \cdot \dfrac{\overline{B_1C}}{\overline{B_1A}} \cdot \dfrac{\overline{C_1A}}{\overline{C_1B}} = -1$ ".

**Theorem 1:** Let us have the polygon $A_1A_2...A_n$, a point $M$ in its plane, and a circular permutation

$$p = \begin{pmatrix} 1 & 2 & ... & n-1 & n \\ 2 & 3 & ... & n & 1 \end{pmatrix}.$$ One notes $M_{ij}$ the intersections of the line $A_iM$ with the lines $A_{i+s}A_{i+s+1},...,A_{i+s+t-1}A_{i+s+t}$ (for all $i$ and $j$, $j \in \{i+s,...,i+s+t-1\}$).

If $M_{ij} \neq A_n$ for all respective indices, and if $2s+t = n$, one has:

$$\prod_{i,j=1, i+s}^{n, i+s+t-1} \frac{\overline{M_{ij}A_j}}{\overline{M_{ij}A_{p(j)}}} = (-1)^n \quad (s \text{ and } t \text{ are natural non-zero numbers}).$$

Analytical proof: Let $M$ be a point in the plain of the triangle $ABC$, such that it satisfies the conditions of the theorem. One chooses a Cartesian system of axes, such that the two parallels with the axes which pass through $M$ do not pass by any point $A_i$ (this is possible).

One considers $M(a,b)$, where $a$ and $b$ are real variables, and $A_i(X_i, Y_i)$ where $X_i$ and $Y_i$ are known, $i \in \{1,2,...,n\}$.

The former choices ensure us the following relations:
$X_i - a \neq 0$ and $Y_i - b \neq 0$ for all $i \in \{1,2,...,n\}$.

The equation of the line $A_iM$ ($1 \le i \le n$) is:
$\dfrac{x-a}{X_i - a} - \dfrac{y-b}{Y_i - b} = 0$. One notes that $d(x,y;X_i,Y_i) = 0$.

One has

$$\frac{\overline{M_{ij}A_j}}{\overline{M_{ij}A_{p(j)}}} = \frac{\delta(A_j, A_iM)}{\delta(A_{p(j)}, A_iM)} = \frac{d(X_j,Y_j;X_i,Y_i)}{d(X_{p(j)},Y_{p(j)};X_i,Y_i)} = \frac{D(j,i)}{D(p(j),i)}$$

where $\delta(A, ST)$ is the distance from $A$ to the line $ST$, and where one notes with $D(a,b)$ for $d(X_a, Y_a; X_b, Y_b)$.

Let's calculate the product, where we will use the following convention: $a+b$ will mean $\underbrace{p(p(...p(a)...))}_{b \text{ times}}$, and $a-b$ will mean $\underbrace{p^{-1}(p^{-1}(...p^{-1}(a)...))}_{b \text{ times}}$

$$\prod_{j=i+s}^{i+s+t-1} \frac{\overline{M_{ij}A_j}}{\overline{M_{ij}A_{j+1}}} = \prod_{j=i+s}^{i+s+t-1} \frac{D(j,i)}{D(j+1,i)} =$$

$$= \frac{D(i+s,i)}{D(i+s+1,i)} \cdot \frac{D(i+s+1,i)}{D(i+s+2,i)} \cdots \frac{D(i+s+t-1,i)}{D(i+s+t,i)} =$$

$$= \frac{D(i+s,i)}{D(i+s+t,i)} = \frac{D(i+s,i)}{D(i-s,i)}$$

The initial product is equal to:
$$\prod_{i=1}^{n} \frac{D(i+s,i)}{D(i-s,i)} = \frac{D(1+s,1)}{D(1-s,1)} \cdot \frac{D(2+s,2)}{D(2-s,2)} \cdots \frac{D(2s,s)}{D(n,s)} \cdot \frac{D(2s+1,s+1)}{D(1,s+1)}$$

$$\cdot \frac{D(2s+2,s+2)}{D(2,s+2)} \cdots \frac{D(2s+t,s+t)}{D(t,s+t)} \cdot \frac{D(2s+t+1,s+t+1)}{D(t+1,s+t+1)} \cdot$$

$$\cdot \frac{D(2s+t+2,s+t+2)}{D(t+2,s+t+2)} \cdots \frac{D(2s+t+s,s+t+s)}{D(t+s,s+t+s)} =$$

$$= \frac{D(1+s,1)}{D(1,1+s)} \cdot \frac{D(2+s,2)}{D(2,2+s)} \cdots \frac{D(2s+t,s+t)}{D(s+t,2s+t)} \cdots \frac{D(s,n)}{D(n,s)} =$$

$$= \prod_{i=1}^{n} \frac{D(i+s,i)}{D(i,i+s)} = \prod_{i=1}^{n} \left( -\frac{P(i+s)}{P(i)} \right) = (-1)^n$$

because:
$$\frac{D(r,p)}{D(p,r)} = \frac{\dfrac{X_r-a}{X_p-a} - \dfrac{Y_r-b}{Y_p-b}}{\dfrac{X_p-a}{X_r-a} - \dfrac{Y_p-b}{Y_r-b}} = -\frac{(X_r-a)(Y_r-b)}{(X_p-a)(Y_p-b)} = -\frac{P(r)}{P(p)},$$

the last equality resulting from what one notes: $(X_t - a)(Y_t - b) = P(t)$. From (1) it results that $P(t) \neq 0$ for all $t$ from $\{1, 2, ..., n\}$. The proof is completed.

**Comments regarding Theorem 1:**

$t$ represents the number of lines of a polygon which are intersected by a line $A_{i_0}M$; if one notes the sides $A_iA_{i+1}$ of the polygon, by $a_i$, then $s+1$ represents the order of the first line intersected by the line $A_1M$ (that is $a_{s+1}$ the first line intersected by $A_1M$).

*Example*: If $s = 5$ and $t = 3$, the theorem says that:
- the line $A_1M$ intersects the sides $A_6A_7, A_7A_8, A_8A_9$.
- the line $A_2M$ intersects the sides $A_7A_8, A_8A_9, A_9A_{10}$.
- the line $A_3M$ intersects the sides $A_8A_9, A_9A_{10}, A_{10}A_{11}$, etc.

*Observation:* The restrictive condition of the theorem is necessary for the existence of the ratios $\dfrac{\overline{M_{ij}A_j}}{\overline{M_{ij}A_{p(j)}}}$.

**Consequence 1.1:** Let's have a polygon $A_1A_2...A_{2k+1}$ and a point $M$ in its plan. For all $i$ from $\{1, 2, ..., 2k+1\}$, one notes $M_i$ the intersection of the line $A_iA_{p(i)}$ with the line which passes through $M$ and by the vertex which is opposed to this line. If $M_i \notin \{A_i, A_{p(i)}\}$ then one has: $\prod_{i=1}^{n} \dfrac{\overline{M_iA_i}}{\overline{M_iA_{p(i)}}} = -1$.

The demonstration results immediately from the theorem, since one has $s = k$ and $t = 1$, that is $n = 2k+1$.

The reciprocal of this consequence is not true.

From where it results immediately that the reciprocal of the theorem is not true either.

Counterexample:

Let us consider a polygon of 5 sides. One plottes the lines $A_1M_3, A_2M_4$ and $A_3M_5$ which intersect in $M$.

Let us have $K = \dfrac{\overline{M_3A_3}}{\overline{M_3A_4}} \cdot \dfrac{\overline{M_4A_4}}{\overline{M_4A_5}} \cdot \dfrac{\overline{M_5A_5}}{\overline{M_5A_1}}$

Then one plots the line $A_4M_1$ such that it does not pass through $M$ and such that it forms the ratio:

(2) $\dfrac{\overline{M_1A_1}}{\overline{M_1A_2}} = 1/K$ or $2/K$. (One chooses one of these values, for which $A_4M_1$ does not pass through $M$).

At the end one traces $A_5M_2$ which forms the ratio $\dfrac{\overline{M_2A_2}}{\overline{M_2A_3}} = -1$ or $-\dfrac{1}{2}$ in function of (2). Therefore the product:

$\prod_{i=1}^{5} \dfrac{\overline{M_iA_i}}{\overline{M_iA_{p(i)}}} = -1$ without having the respective lines concurrent.

**Consequence 1.2:** Under the conditions of the theorem, if for all i and $j, j \notin \{i, p^{-1}(i)\}$, one notes $M_{ij} = A_i M \cap A_j A_{p(j)}$ and $M_{ij} \notin \{A_j, A_{p(j)}\}$ then one has:

$$\prod_{i,j=1}^{n} \frac{\overline{M_{ij} A_j}}{\overline{M_{ij} A_{p(j)}}} = (-1)^n .$$

$j \notin \{i, p^{-1}(i)\}$
Effectively one has $s = 1$, $t = n - 2$, and therefore $2s + t = n$.

**Consequence 1.3:** For $n = 3$, it comes $s = 1$ and $t = 1$, therefore one obtains (as a particular case) the theorem of Ceva.

**An Application of the Generalizations of Ceva's Theorem** is presented below.

**Theorem 2:** Let us consider a polygon $A_1 A_2 ... A_n$ inserted in a circle. Let $s$ and $t$ be two non zero natural numbers such that $2s + t = n$. By each vertex $A_i$ passes a line $d_i$ which intersects the lines $A_{i+s} A_{i+s+1}, ...., A_{i+s+t-1} A_{i+s+t}$ at the points $M_{i,i+s}, ..., M_{i+s+t-1}$ respectively and the circle at the point $M_i'$. Then one has:

$$\prod_{i=1}^{n} \prod_{j=i+s}^{i+s+t-1} \frac{\overline{M_{ij} A_j}}{\overline{M_{ij} A_{j+1}}} = \prod_{i=1}^{n} \frac{\overline{M_i' A_{i+s}}}{\overline{M_i' A_{i+s+t}}} .$$

*Proof:*
Let $i$ be fixed.
1) The case where the point $M_{i,i+s}$ is inside the circle.

There are triangles $A_i M_{i,i+s} A_{i+s}$ and $M_i' M_{i,i+s} A_{i+s+1}$ which are similar, since the angles $M_{i,i+s} A_i A_{i+s}$ and $M_{i,i+s} A_{i+s+1} M_i'$ on one side, and $A_i M_{i,i+s} A_{i+s}$ and $A_{i+s+1} M_{i,i+s} M_i'$ are equal. It results from it that:

(1) $\quad \dfrac{\overline{M_{i,i+s} A_i}}{\overline{M_{i,i+s} A_{i+s+1}}} = \dfrac{\overline{A_i A_{i+s}}}{\overline{M_i' A_{i+s+1}}}$

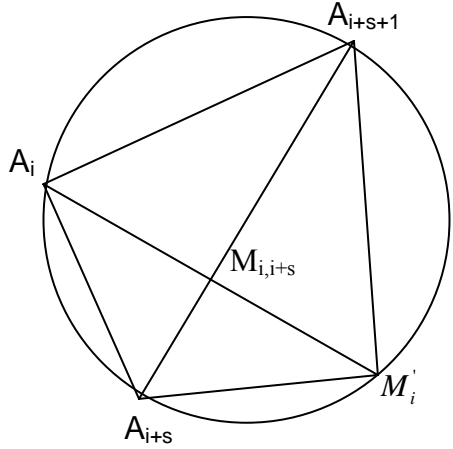

In a similar manner, one shows that the triangles $M_{i,i+s}A_iA_{i+s+1}$ and $M_{i,i+s}A_{i+s}M'_i$ are similar, from which:

(2) $\quad \dfrac{\overline{M_{i,i+s}A_i}}{\overline{M_{i,i+s}A_{i+s}}} = \dfrac{\overline{A_iA_{i+s+1}}}{\overline{M'_iA_{i+s}}}$. Dividing (1) by (2) we obtain:

(3) $\quad \dfrac{\overline{M_{i,i+s}A_{i+s}}}{\overline{M_{i,i+s}A_{i+s+1}}} = \dfrac{\overline{M'_iA_{i+s}}}{\overline{M'_iA_{i+s+1}}} \cdot \dfrac{\overline{A_iA_{i+s}}}{\overline{A_iA_{i+s+1}}}$.

2) The case where $M_{i,i+s}$ is exterior to the circle is similar to the first, because the triangles (notations as in 1) are similar also in this new case. There are the same interpretations and the same ratios; therefore one has also the relation (3).

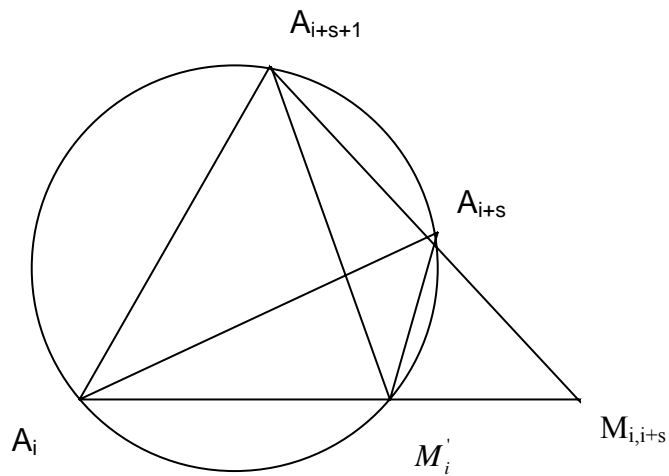

Let us calculate the product:

$$\prod_{j=i+s}^{i+s+t-1} \frac{\overline{M_{ij}A_j}}{\overline{M_{ij}A_{j+1}}} = \prod_{j=i+s}^{i+s+t-1} \left( \frac{\overline{M_i'A_j}}{\overline{M_i'A_{j+1}}} \cdot \frac{\overline{A_iA_j}}{\overline{A_iA_{j+1}}} \right) =$$

$$= \frac{\overline{M_i'A_{i+s}}}{\overline{M_i'A_{i+s+1}}} \cdot \frac{\overline{M_i'A_{i+s+1}}}{\overline{M_i'A_{i+s+2}}} \cdots \frac{\overline{M_i'A_{i+s+t-1}}}{\overline{M_i'A_{i+s+t}}} \cdot$$

$$\cdot \frac{\overline{A_iA_{i+s}}}{\overline{A_iA_{i+s+1}}} \cdot \frac{\overline{A_iA_{i+s+1}}}{\overline{A_iA_{i+s+2}}} \cdots \frac{\overline{A_iA_{i+s+t-1}}}{\overline{A_iA_{i+s+t}}} = \frac{\overline{M_i'A_{i+s}}}{\overline{M_i'A_{i+s+t}}} \cdot \frac{\overline{A_iA_{i+s}}}{\overline{A_iA_{i+s+t}}}$$

Therefore the initial product is equal to:

$$\prod_{i=1}^{n} \left( \frac{\overline{M_i'A_{i+s}}}{\overline{M_i'A_{i+s+t}}} \cdot \frac{\overline{A_iA_{i+s}}}{\overline{A_iA_{i+s+t}}} \right) = \prod_{i=1}^{n} \frac{\overline{M_i'A_{i+s}}}{\overline{M_i'A_{i+s+t}}}$$

since:

$$\prod_{i=1}^{n} \frac{\overline{A_iA_{i+s}}}{\overline{A_iA_{i+s+t}}} = \frac{\overline{A_1A_{1+s}}}{\overline{A_1A_{1+s+t}}} \cdot \frac{\overline{A_2A_{2+s}}}{\overline{A_2A_{2+s+t}}} \cdots \frac{\overline{A_sA_{2s}}}{\overline{A_{s+1}A_1}} \cdot$$

$$\cdot \frac{\overline{A_{s+2}A_{2s+2}}}{\overline{A_{s+2}A_2}} \cdots \frac{\overline{A_{s+t}A_n}}{\overline{A_{s+t}A_t}} \cdot \frac{\overline{A_{s+t+1}A_1}}{\overline{A_{s+t+1}A_{t+1}}} \cdot \frac{\overline{A_{s+t+2}A_2}}{\overline{A_{s+t+2}A_{t+2}}} \cdots \frac{\overline{A_nA_s}}{\overline{A_nA_{s+t}}} = 1$$

(by taking into account the fact that $2s+t=n$).

**Consequence 2.1:** If there is a polygon $A_1A_2,...,A_{2s-1}$ inscribed in a circle, and from each vertex $A_i$ one traces a line $d_i$ which intersects the opposite side $A_{i+s-1}A_{i+s}$ in $M_i$ and the circle in $M_i'$ then:

$$\prod_{i=1}^{n} \frac{\overline{M_iA_{i+s-1}}}{\overline{M_iA_{i+s}}} = \prod_{i=1}^{n} \frac{\overline{M_i'A_{i+s-1}}}{\overline{M_i'A_{i+s}}}$$

In fact for $t=1$, one has $n$ odd and $s = \frac{n+1}{2}$.

If one makes $s=1$ in this consequence, one finds the mathematical note from [1], pages 35-37.

**Application:** If in the theorem, the lines $d_i$ are concurrent, one obtains:

$$\prod_{i=1}^{n} \frac{\overline{M_i'A_{i+s}}}{\overline{M_i'A_{i+s+t}}} = (-1)^n.$$

**Reference:**

[1]   Dan Barbilian - Ion Barbu – "Pagini inedite", Editura Albatros, Bucharest, 1981 (Ediție îngrijită de Gerda Barbilian, V. Protopopescu, Viorel Gh. Vodă).